%% file: main.tex
\documentclass[letterpaper, 10 pt, conference]{ieeeconf}  

\IEEEoverridecommandlockouts                              

\usepackage{graphicx}
\usepackage{subcaption}
\usepackage{latexsym,dsfont}
\usepackage{amsmath,amssymb,amsfonts}
\usepackage{xcolor}
\usepackage{booktabs}
\usepackage{hyperref}

\usepackage{tikz}
\input{abbrv.tex}

\graphicspath{{figures/}{./}}

\newtheorem{assumption}{Standing Assumption}

\title{\LARGE \bf
A Neural Network Approach Applied to Multi-Agent Optimal Control
}

\author{Derek Onken, Levon Nurbekyan, Xingjian Li, Samy Wu Fung, Stanley Osher, and Lars Ruthotto%
\thanks{This work was supported in part by NSF award DMS 1751636, AFOSR Grants 20RT0237 \& FA9550-18-1-0167, AFOSR MURI FA9550-18-1-0502,
Binational Science Foundation Grant 2018209,
US DOE, Office of Advanced Scientific Computing Research Field Work Proposal 20-023231, ONR Grants No.
N00014-18-1-2527 \& N00014-20-1-2093, a gift from UnitedHealth Group R\&D, and a GPU donation by NVIDIA Corporation.}%
\thanks{D.~Onken is with the Department of Computer Science,
Emory University, Atlanta, GA, USA (email: {\tt\small donken@emory.edu})}%
\thanks{L.~Nurbekyan, S.~Wu Fung, and S.~Osher are with the Department of Mathematics,
UCLA, Los Angeles, CA, USA (email: {\tt\small \{lnurbek; swufung; sjo\}@math.ucla.edu})}%
\thanks{X.~Li and L.~Ruthotto are with the Department of Mathematics,
Emory University, Atlanta, GA, USA (email: {\tt\small \{xingjian.li; lruthotto\}@emory.edu})}%
}

\begin{document}

\maketitle
\thispagestyle{empty}
\pagestyle{empty}

\begin{abstract}

We propose a neural network approach for solving high-dimensional optimal control problems.
In particular, we focus on multi-agent control problems with obstacle and collision avoidance.
These problems immediately become high-dimensional, even for moderate phase-space dimensions per agent.
Our approach fuses the Pontryagin Maximum Principle and Hamilton-Jacobi-Bellman (HJB) approaches and
parameterizes the value function with a neural network. Our approach yields controls in a feedback form
for quick calculation and robustness to moderate disturbances to the system.
We train our model using the objective function and optimality conditions of the control problem. 
Therefore, our training algorithm neither involves a data generation phase nor solutions from another
algorithm. Our model uses empirically effective HJB penalizers for efficient training.
By training on a distribution of
initial states, we ensure the controls' optimality is achieved on a large portion of the state-space.
Our approach is grid-free and scales efficiently to dimensions where grids become impractical or infeasible.
We demonstrate our approach's effectiveness on a 150-dimensional multi-agent problem with obstacles.

\end{abstract}

\section{Introduction} \label{sec:introduction}

Optimal control (OC) problems are ubiquitous in pure and applied mathematics, physics, computer science, engineering, finance, and elsewhere \cite{flemingsoner06,bardi97,cannarsa04}. Thus, theoretical and numerical analyses of OC problems have paramount importance across disciplines. We focus on developing numerical solution methods for general high-dimensional OC problems. 

Two of the most common strategies to solve OC problems are Pontryagin's Maximum Principle (PMP)~\cite{pontryagin62} and the Hamilton-Jacobi-Bellman (HJB) PDE~\cite{bellman57} (Sec.~\ref{sec:preliminaries}).

PMP is a \textit{local} solution method because optimal controls correspond to fixed initial states.
Necessary conditions render an ODE system for the state and \textit{adjoint} variables and a \textit{maximum principle} relating the adjoint variable with the optimal control. 
This approach is grid-free and thus suitable for high-dimensional problems. However, the aforementioned ODE system is often challenging to solve due to its forward-backward structure \cite{kang2017mitigating,kang2019algorithms}.
Additionally, the OC problem's non-convexity renders the possibility of multiple non-optimal solutions, and additional considerations are necessary, such as the differentiability of the value function \cite[Theorem 7.3.9]{cannarsa04}. Unfortunately, these conditions are virtually impossible to enforce or verify numerically. The effects of non-convexity are especially pronounced in multi-agent collision-avoidance problems \cite[I.A]{riviere2020glas}. For OC problems that are convex, high-dimensional solvers can be devised via primal-dual convex optimization methods \cite{lin2018splitting,darbonosher16,Kirchner18,Kirchner2018APM}.

Because PMP is a local solution method,
shocks or sudden changes in the system's initial conditions lead to a new optimization problem. PMP is unattractive for real-time applications because the control-search time is vital.

Alternatively, the HJB approach aims at solving the OC problem for \textit{all initial states at once}. More precisely, the value function---also known as the optimal cost-to-go---of an OC problem is a solution of a suitable HJB PDE.
After computing the value function, we can recover the optimal control at any state from the value function's gradient. Such controls are said to be in feedback form and are especially useful for real-time applications if calculation-times for the feedback form are short.

Although effective, HJB equations are challenging to solve numerically, especially when the state's dimension $d \geq 4$. First, in a deterministic setup, HJB is a first-order non-linear equation and generally does not admit smooth solutions. Second, traditional numerical methods for HJB equations, such as ENO/WENO \cite{osher1991high}, rely on grids and therefore suffer from the curse of dimensionality \cite{bellman57}.

We propose a machine learning framework to overcome the curse of dimensionality by approximating the value function with a neural network (NN). Our approach is a fusion of PMP and HJB.
Combining the PMP and HJB approaches, we express the control cost in terms of the value function and search for an NN approximation that minimizes this cost on a cloud of initial states.
Furthermore, we improve the NN training by adding HJB residual penalties (Fig.~\ref{fig:penalizers}), similar to \cite{ruthotto2020machine,lin2020apac,onken2020otflow}. 

Our approach has several advantages. First, it applies to generic OC problems. Second, we find controls in a feedback form that is crucial for real-time applications. Training the NN on a cloud of initial states ensures the controls' optimality on a large portion of the state-space.
Consequently, controls are robust to moderate disturbances or shocks to the system (Fig.~\ref{fig:shock}). Third, our method is grid-free and suitable for high-dimensional problems. Finally, our approach incorporates machine learning techniques for solving large scale optimization problems. 
As a result, we are able to solve $150$-dimensional OC problems (Sec.~\ref{subsec:swarm}).

Applications of deep learning techniques to OC problems appear in seminal works \cite{han2016deep,E_2017,Han_2018}, where the authors apply backward stochastic differential equations techniques to solve high-dimensional stochastic OC problems. In \cite{nusken2020solving}, the authors extend these methods by introducing and analyzing different loss functions. These works consider stochastic problems with fixed initial states. We focus on finding solutions to deterministic problems that are robust to shocks.

In \cite{nakamura2019adaptive}, the authors first generate optimal controls for a sample of initial states then train an NN to fit this data. The data generation phase is performed by a different algorithm \cite{kang2019algorithms}. A similar approach exists in \cite{riviere2020glas} with a different data generation algorithm \cite{honig2018trajectory}. In contrast, our approach directly minimizes the cost function without a generation phase.

Our work is based on the same framework as \cite{kunisch2020semiglobal}, which approximates the feedback control with an NN then optimizes the control cost on a cloud of initial states and provides a theoretical analysis of OC solutions via NN approximations. We extend the framework to finite horizon problems with non-quadratic costs and parameterize the value function instead of the feedback function. This latter approach enables enforcing HJB conditions, which empirically improves numerical performance for solving high-dimensional mean-field games, mean-field control, and normalizing flows~\cite{ruthotto2020machine,lin2020apac,onken2020otflow}. We demonstrate similar advantages in the OC problems considered in this work.

We consider deterministic OC problems with a particular emphasis on centrally controlled multi-agent systems. For $n$ agents in an $q$-dimensional space we obtain a $d=n\cdot q$-dimensional OC problem. Thus, even moderate $n,q$ yield problems that are out of reach for traditional HJB solvers. 

We demonstrate the effectiveness of our method by solving a $50$-agent control problem in a $3$-dimensional space with obstacle and collision avoidance (Fig.~\ref{fig:quadtraj}). 
The overall dimension of this problem is $d{=}150$. 
Additionally, we demonstrate  our model's robustness to shocks and the effect of the penalizers on a $4$-dimensional corridor problem.

\section{Preliminaries} \label{sec:preliminaries}

Here, we briefly recall relevant OC theory. We refer to \cite[Chapters I, II]{flemingsoner06} for a detailed exposition. We are interested in deterministic fixed finite time-horizon problems. Consider the time-horizon $[0,T]$ and the system's dynamics given by
\begin{equation}\label{eq:charOC}
    \partial_s \bfz(s) = f(s,\bfz(s),\bfu(s)), \quad t \leq  s \leq  T, \quad \bfz(t)= \bfx.
\end{equation}
Here, $\bfz \in \R^d$ describes the state of the system, and $\bfu \in U \subset \R^a$ describes the controls. Hence, $f \colon [0,T]\times \R^d\times U \to \R^d$ models the evolution of the state $\bfz \colon [t,T]\to \R^d$ in response to the application of a control $\bfu \colon [t,T]\to U$ for initial time $t\in [0,T]$ and initial state $\bfx$. We assume that $f,L,G,U$ are sufficiently regular (see \cite[Sec. I.3, I.8-9]{flemingsoner06} for a list of assumptions). Next, assume that the control $\bfu$ yields cost
\begin{equation}\label{eq:Joc}
    J_{t,\bfx} [\bfu] = G\big(\bfz(T)\big) + \int_{t}^T L\big(s,\bfz(s), \bfu(s)\big) \, \du s,
\end{equation}
where $L
\colon [0,T]\times \R^d \times U \to \R$ is the \textit{running cost} or the \textit{Lagrangian}, and $G \colon \R^d \to \R$ is the \textit{terminal cost}. OC problems seek the control that incurs the minimal cost; i.e.,
\begin{equation}\label{eq:OC}
    \Phi(t,\bfx)=\inf_{\bfu} J_{t,\bfx}[\bfu],
\end{equation}
where $\Phi$ is called the \textit{value function}. A solution $\bfu^*$ of \eqref{eq:OC} is called an \textit{optimal control}. Accordingly, the $\bfz^*$ which corresponds to $\bfu^*$ is called an \textit{optimal trajectory}.

Next, the \textit{Hamiltonian} of the system is given by
\begin{equation}\label{eq:H}
\begin{split}
    H(t,\bfx,\bfp)=&\sup_{\bfu \in U} \left\{ -\bfp \cdot f(t,\bfx,\bfu)-L(t,\bfx,\bfu) \right\},
\end{split}
\end{equation}
where $\bfp$ is the \textit{adjoint state}. The following is a standing assumption throughout the paper.
\begin{assumption}\label{assum:feedback}
Assume that \eqref{eq:H} admits a unique continuous closed-form solution $\bfu^*(t,\bfx,\bfp)$.
\end{assumption}
Under this assumption and denoting $\nabla$ as the gradient with respect to the state variables only, one has that%
\begin{equation}\label{eq:L(u*)}%
    L(t,\bfx,\bfu^*(t,\bfx,\bfp))=\bfp\cdot \nabla H(t,\bfx,\bfp)-H(t,\bfx,\bfp).
\end{equation}

The PMP \cite{pontryagin62} states that for a solution $(\bfz^*,\bfu^*)$ of \eqref{eq:OC} there exist $\bfp \colon [0,T]\to \R^d$ such that, for $t\leq s \leq T$, 
\begin{equation}\label{eq:PMP_combined}
    \begin{cases}
    \partial_s \bfz^*(s)=-\nabla_{\bfp} H \big(s,\bfz^*(s),\bfp(s) \big),\\
    \partial_s \bfp(s)=\nabla H \big( s,\bfz^*(s),\bfp(s) \big),\\
    \bfz^*(t)=\bfx,\quad \bfp(T)=\nabla G \big( \bfz^*(T) \big),\\
    \bfu^*(s) =\bfu^* \big(s,\bfz^*(s),\bfp(s) \big).
    \end{cases}
\end{equation}

The HJB PDE, also known as the \textit{dynamic programming} equation, corresponding to \eqref{eq:OC} is given by 
\begin{equation}\label{eq:HJB}
-\partial_t \Phi(t,\bfx) + H\big(t,\bfx,\nabla \Phi(t,\bfx)\big)=0, ~ \Phi(T,\bfx) = G(\bfx),
\end{equation}
for $(t,\bfx)\in [0,T] \times \R^d $. The cornerstone of the HJB approach is that the value function $\Phi$ is the unique viscosity solution of \eqref{eq:HJB}. Moreover, $\bfp$ in \eqref{eq:PMP_combined} and $\Phi$ are related by
\begin{equation}\label{eq:gradPhi}
\bfp(s)=\nabla \Phi \big(s,\bfz^*(s) \big),\quad t<s\leq T.
\end{equation}
Thus, the optimal control $\bfu^*$ is given in a feedback form
\begin{equation}\label{eq:opt_control_nec}
    \bfu^*(s)=\bfu^*\left(s,\bfz^*(s),\nabla \Phi \big(s,\bfz^*(s) \big)\right),\quad t<s\leq T,
\end{equation}
and $\bfz^*$ evolves according to
\begin{equation}\label{eq:z_evol}
\begin{cases}
\partial_s \bfz^*(s)=-\nabla_{\bfp} H\left(s,\bfz^*(s),\nabla \Phi \big(s,\bfz^*(s) \big)\right), \\
\bfz^*(t)=\bfx.
\end{cases} 
\end{equation}

\section{Machine Learning Approach} \label{sec:formulation}

We derive a machine learning formulation of \eqref{eq:OC} by leveraging the advantages of both the PMP and the HJB approaches. In particular, we directly optimize~\eqref{eq:Joc} subject to~\eqref{eq:charOC}. 
Rather than solving for the controls, we first postulate the dependence of $(\bfz,\bfu)$ on $\Phi$ according to~\eqref{eq:opt_control_nec},~\eqref{eq:z_evol} and parameterize $\Phi$ by an NN. We also add penalizers that punish deviations from \eqref{eq:HJB}, similar to \cite{ruthotto2020machine,onken2020otflow}.

\subsection{Main formulation}

Denote the control variable corresponding to initial position $\bfx \in \R^d$ at time $t{=}0$ by $\bfu_{\bfx} \colon [0,T]\to U$. And denote the trajectory corresponding to initial data $(0,\bfx)$ by $\bfz_{\bfx}$; i.e.,
\begin{equation*}
    \partial_s \bfz_{\bfx}(s) = f(s,\bfz_{\bfx}(s),\bfu_{\bfx}(s)), \quad 0 \leq  s \leq  T, \quad \bfz_{\bfx}(0)= \bfx.
\end{equation*}
Furthermore, fix a $\rho_0 \in \CP(\R^d)$ and consider the problem
\begin{equation}\label{eq:cont_loss}
\begin{split}
\inf_{\{\bfu_{\bfx}\}} \int_{\R^d} J_{0,\bfx} [\bfu_{\bfx}] \rho_0(\bfx) \, \du \bfx   = \inf_{\{\bfu_{\bfx}\}} \E_{\bfx \sim \rho_0} J_{0,\bfx} [\bfu_{\bfx}].
\end{split}
\end{equation}
Note that $\{\bfu^*_{\bfx}\}$ solves \eqref{eq:cont_loss} if and only if $\bfu^*_{\bfx}$ is an optimal control for $\rho_0$ almost everywhere $\bfx \in \R^d$. Thus, solving \eqref{eq:OC} with initial points $\bfx$ in some domain $\Omega$ is equivalent to solving \eqref{eq:cont_loss} for $\rho_0 \in \CP(\R^d)$ such that $\operatorname{supp}(\rho_0){=}\Omega$. 
Formulation \eqref{eq:cont_loss} is employed in mean-field game and control systems~\cite{ruthotto2020machine,lin2020apac} and in stabilization problems~\cite{kunisch2020semiglobal}.

Postulating \eqref{eq:opt_control_nec}, \eqref{eq:z_evol}, we arrive at our main formulation
\begin{equation}\label{eq:cont_loss_Phi}
\begin{split}
    \inf_{\Phi} &~\E_{\bfx \sim \rho_0} \bigg( G\big(\bfz_{\bfx}(T)\big) + \int_{0}^T L \big(s,\bfz_{\bfx}(s),\bfu_{\bfx}(s) \big)\, \du s\bigg)\\
\mbox{s.t.} &~
\begin{cases}
\bfu_{\bfx}(s)=\bfu^*\left(s,\bfz_{\bfx}(s),\nabla \Phi \big(s,\bfz_{\bfx}(s) \big)\right),\\
\partial_s \bfz_{\bfx}(s)=-\nabla_{\bfp} H\left(s,\bfz_{\bfx}(s),\nabla \Phi \big(s,\bfz_{\bfx}(s) \big)\right),\\
\bfz_{\bfx}(0)=\bfx.
\end{cases}
\end{split}
\raisetag{30pt}
\end{equation}

In each OC problem, we set up and solve~\eqref{eq:cont_loss_Phi}, where $f$, $L$, and thus $H$ vary with the problem. The $L$ contains terms that reflect the features of the problem (Sec.~\ref{sec:lagrangian}).

We approximate $\Phi(\cdot)$ by an NN, $\Phi(\cdot;\bfth)$ (Sec.~\ref{sec:model}), and turn \eqref{eq:cont_loss_Phi} into a finite-dimensional optimization over the weights $\bfth$. For brevity, we often omit $\Phi$'s explicit dependence on $\bfth$.

\subsection{Adding HJB penalizers}

We introduce three penalty terms $c_{\rm HJt,\bfx}$, $c_{\rm HJfin,\bfx}$, and $c_{\rm HJgrad,\bfx}$ derived from the HJB PDE~\eqref{eq:HJB} as follows:
\begin{equation} \label{eq:HJt}
\begin{split}
    c_{\rm HJt,\bfx}(t) =& \int_0^t  \big|\partial_s \Phi \big( s,\bfz_{\bfx}(s) ; \bfth \big) \\
    &- H \big( s,\bfz_{\bfx}(s),\nabla \Phi(s,\bfz_{\bfx}(s); \bfth ) \big) \big| \du s \\
    c_{\rm HJfin,\bfx} =& |\Phi(T,\bfz_{\bfx}(T) ; \bfth ) - G(\bfz_{\bfx}(T)) | \\
    c_{\rm HJgrad,\bfx} =& |\nabla \Phi(T, \bfz_{\bfx}(T) ; \bfth ) - \nabla G(\bfz_{\bfx}(T))| .
\end{split}
\end{equation}
The ${\rm HJ_{t}}$ penalizer arises from the first equation in~\eqref{eq:HJB}, whereas ${\rm HJ_{fin}}$ and ${\rm HJ_{grad}}$ are direct results of the final-time condition in~\eqref{eq:HJB} and its gradient, respectively.
Penalizers prove helpful in training NNs for solving problems similar to~\eqref{eq:cont_loss_Phi}~\cite{ruthotto2020machine,lin2020apac,onken2020otflow,finlay2020train}. 
They improve the training convergence (Sec.~\ref{sec:penalizers}) without altering the solution of \eqref{eq:cont_loss_Phi}.

Adding $c_{\rm HJt,\bfx}, c_{\rm HJfin,\bfx}, c_{\rm HJgrad,\bfx}$ to \eqref{eq:cont_loss_Phi} and rewriting the time-integral in terms of ODE constraints, we obtain
\begin{equation} \label{eq:full_opt}
		\begin{split}
		\min_{\bfth} \; \E_{\bfx \sim \rho_0} \;\; \big( \ell_{\bfx}(T)   + G(\bfz_{\bfx}(T))
		+ \beta_1 c_{\rm HJt,\bfx}(T) \\ \quad + \beta_2 c_{\rm HJfin,\bfx} + \beta_3 c_{\rm HJgrad,\bfx} \big),
		\end{split}
	\end{equation}
	subject to
	\begin{equation}\label{eq:ODE_cons}
		\renewcommand*{\arraystretch}{1.1} 
		\begin{split}
		\partial_s\begin{pmatrix}
		\bfz_{\bfx}(s)\\
		\ell_{\bfx}(s)\\
		c_{\rm HJt,\bfx}(s) \\
		\end{pmatrix} 
		=\begin{pmatrix}
		- \nabla_{\bfp} H(s,\bfz_{\bfx}(s),\nabla \Phi(s, \bfz_{\bfx}(s); \bfth ))\\
		 L_{\bfx}(s)\\
		R_{\bfx}(s)
		\end{pmatrix},
		\end{split}
	\end{equation}
initialized with $\bfz_{\bfx}(0) = \bfx$ and $\ell_{\bfx}(0)= c_{\rm HJt,\bfx}(0) = 0$, and 
\begin{equation*}
\begin{split}
    L_{\bfx}(s)=&\nabla \Phi \big( s, \bfz_{\bfx}(s) ; \bfth \big) \cdot \nabla_{\bfp} H \big( s,\bfz_{\bfx}(s),\nabla \Phi(s, \bfz_{\bfx}(s) ; \bfth) \big)\\
    & -H \big( s,\bfz_{\bfx}(s),\nabla \Phi(s,\bfz_{\bfx}(s) ; \bfth) \big)\\
    R_{\bfx}(s)=&|\partial_s \Phi(s,\bfz_{\bfx}(s) ; \bfth) - H \big( s,\bfz_{\bfx}(s),\nabla \Phi(s, \bfz_{\bfx}(s) ; \bfth) \big)|.
\end{split}
\end{equation*}
We note that reformulating the Lagrangian $L_{\bfx}(s)$ uses \eqref{eq:L(u*)}.

The objective function thus contains the accumulated running cost $\ell_{\bfx}(T)$, the HJB penalty along the trajectories $c_{\rm HJt,\bfx}(T)$, the final-time HJB penalty $c_{\rm HJfin,\bfx}$, and the transversality penalty $c_{\rm HJgrad,\bfx}$. The penalty multipliers $\beta_1,\beta_2,\beta_3>0$ are hyperparameters of the model (Sec.~\ref{sec:model}).

\begin{figure}
	\centering
	\begin{subfigure}{\linewidth}
	  \centering
	  \includegraphics[width=0.93\textwidth]{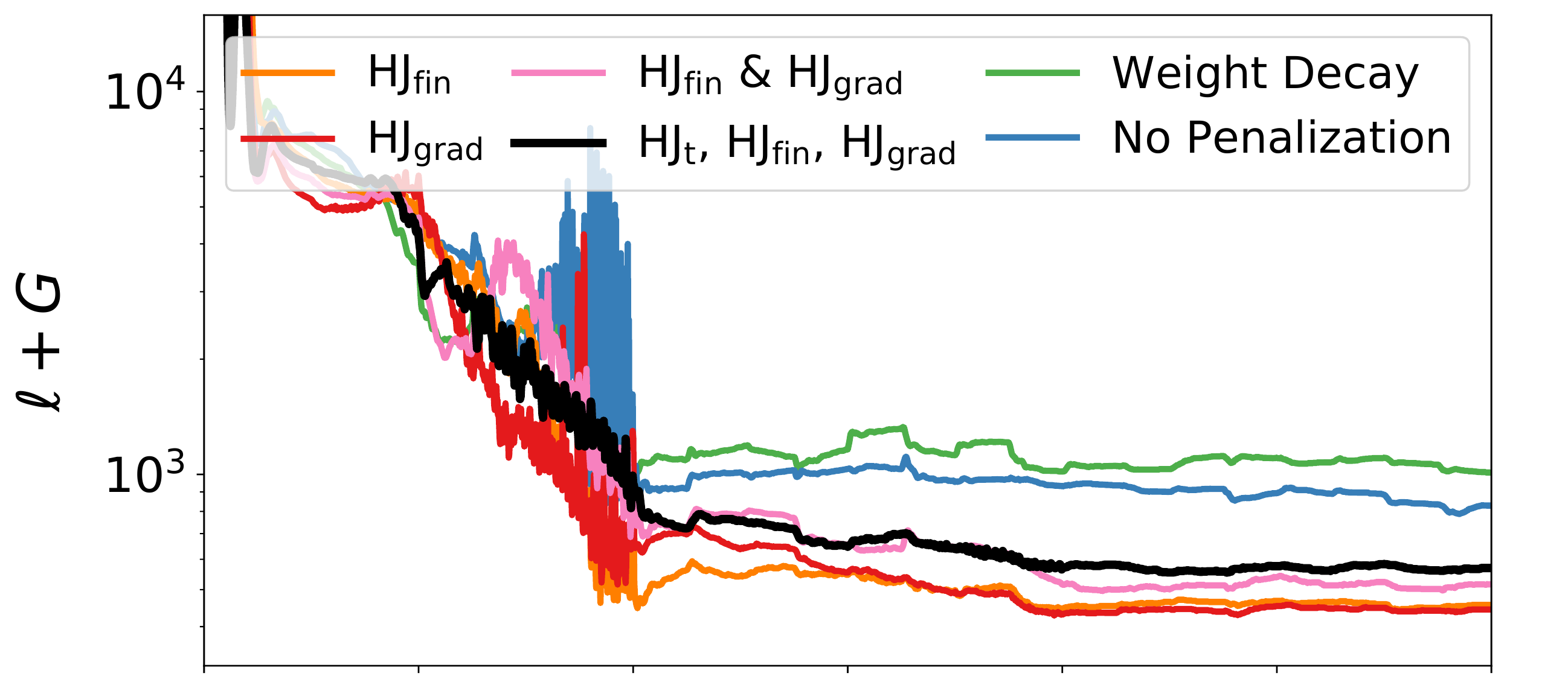}
	\end{subfigure} \\
	\begin{subfigure}{\linewidth}
	  \centering 
	  \includegraphics[width=0.93\textwidth]{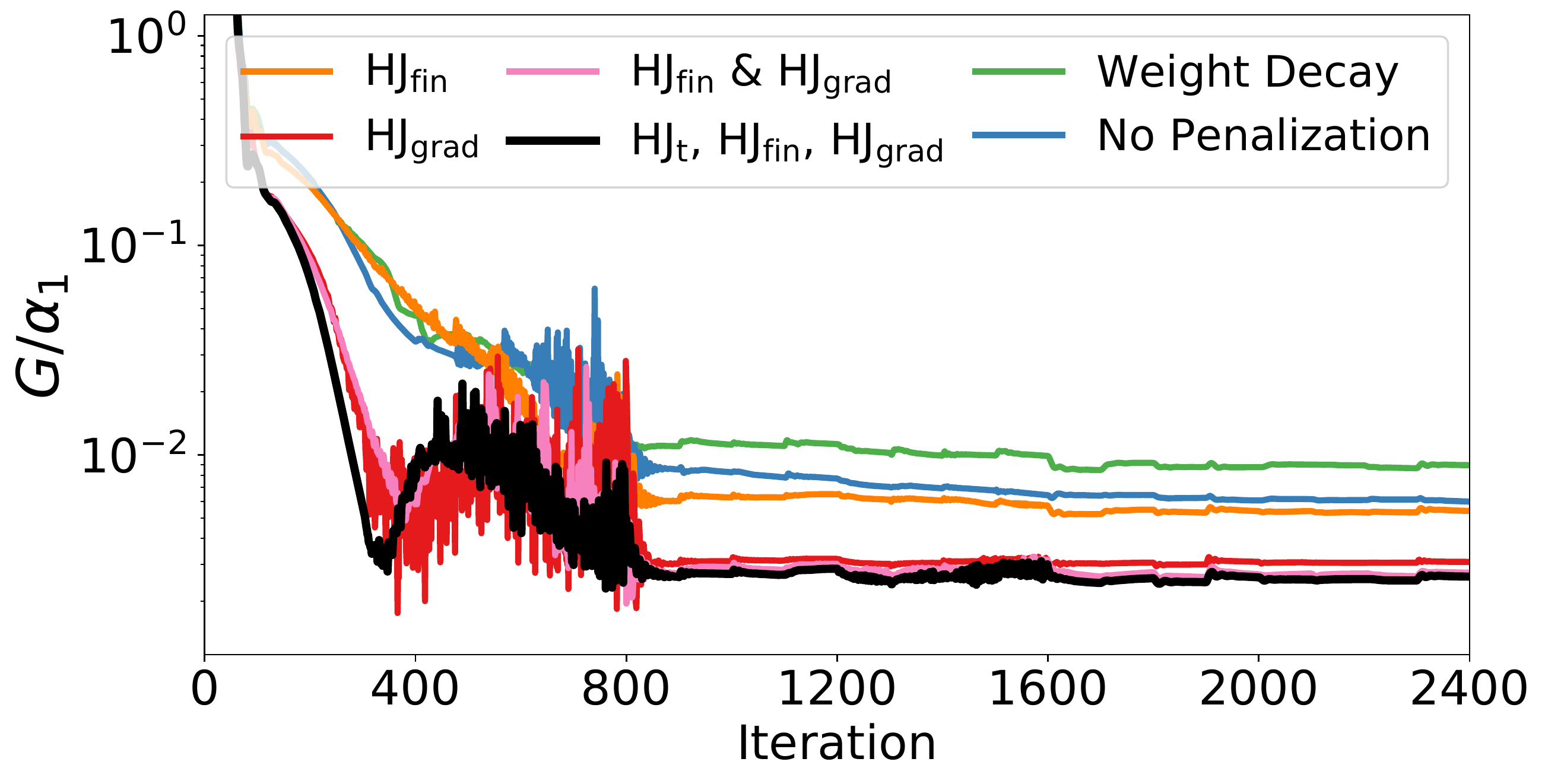}
	\end{subfigure}
	\caption{We compare the validation values for models trained with different penalizers. Using the HJB penalizers leads to quick convergence and a low $G$ value. Each curve is the average of three training instances.}
	\label{fig:penalizers}
\end{figure}

\subsection{The effect of the HJB penalizers} \label{sec:penalizers}

We experimentally assess the effectiveness of the penalizers $c_{\rm HJt}$, $c_{\rm HJfin}$, $c_{\rm HJgrad}$. To this end, we define six models (various combinations of the three HJB penalizers and one with weight decay) and train each on the corridor problem (Sec.~\ref{sec:softcorridor}). Using the HJB penalizers results in the quicker model convergence on a hold-out validation set (Fig.~\ref{fig:penalizers}).

$\mathbf{HJ_{\rm t}}\colon$
We enforce the PDE~\eqref{eq:HJB} describing the time derivative of $\Phi$ along the trajectories. Including this penalizer improves regularity and reduces the necessary number of time steps when solving the dynamics~\cite{ruthotto2020machine,lin2020apac,onken2020otflow,yang2019}.

$\mathbf{HJ_{\rm fin}}\colon$ We enforce the final-time condition of the PDE~\eqref{eq:HJB}. The inclusion of this penalizer helps the network achieve the target~\cite{ruthotto2020machine}. Experimentally, using ${\rm HJ_{fin}}$ correlates with a slightly lower $G$ value (Fig.~\ref{fig:penalizers}). 

$\mathbf{HJ_{\rm grad}}\colon$ We enforce the tranversality condition $\nabla \Phi(T,\bfz(T)){=}\nabla G(\bfz(T)),~\forall \bfz$, a consequence of the final-time HJB condition~\eqref{eq:HJB}. Numerically, all conditions are enforced on a finite sample set. Therefore, higher-order regularization may help the generalization; i.e., achieving a better match of $\Phi(T,\cdot)$ and $G$ for samples not used during training (the hold-out validation set). We observe the latter experimentally; using ${\rm HJ_{grad}}$ instead of ${\rm HJ_{fin}}$ results in lower $G$ (Fig.~\ref{fig:penalizers}). ~\cite{nakamura2019adaptive} similarly enforces $\nabla \Phi$ values.

\subsection{Lagrangian for Obstacle and Collision Avoidance} \label{sec:lagrangian}

For multi-agent problems, the Lagrangian consists of three terms: an energy term $E\colon  U \to \R$ penalizing how much the agents travel, an obstacle term $Q\colon \R^d \to \R$ penalizing agents at certain spatial locations (i.e., a terrain function), and an interaction term $W\colon \R^d \to \R$ penalizing the proximity among agents (i.e, collision avoidance). 

The $n$ agents have initial states $x_1,\dots,x_n \in \R^{q}$.
For initial joint-state $\bfx = (x_1,\dots,x_n) \in \R^d$ with $d=q\cdot n$, we represent the joint-state $\bfz_{\bfx} \in \R^d$ at time $t$ as
\begin{equation}
    \bfz_{\bfx}(t) = (z_{x_1}(t),z_{x_2}(t),\dots,z_{x_n}(t)) ,
\end{equation}
where $z_{x_i} \in \R^{q}$ is the $i$th agent's state. The control follows
\begin{equation}
    \bfu_{\bfx}(t) = (u_{x_1}(t),u_{x_2}(t),\dots,u_{x_n}(t)).
\end{equation}

As a result, we define the Lagrangian as
\begin{equation}\label{eq:mult_ag_obj}
    \begin{split}
        L \big( t, \bfz_{\bfx}, \bfu_{\bfx} \big) = \, E \big( \bfu_{\bfx} \big) + \, \alpha_2 Q \big( \bfz_{\bfx} \big) 
        + \, \alpha_3 W \big( \bfz_{\bfx} \big) \qquad \quad
        \\
        =  \sum_{i=1}^n  E_i \big( u_{x_i} \big) \, + \, \alpha_2 \sum_{i=1}^n  Q_i \big( z_{x_i} \big) + \, \alpha_3 \sum_{j\neq i} W_{ij}\big(z_{x_i},z_{x_j}\big),
    \end{split}
    \raisetag{45pt} 
\end{equation} 
where we omit the dependence on $t$ for brevity.
The scalars $\alpha_2,\alpha_3$ calibrate the magnitude (relative to $E$) of the penalization of the obstacle and interactions, respectively.

\section{Implementation} \label{sec:model}

	We parameterize the value function as
	\begin{equation}
	\label{eq:NNArchitecture}
		\begin{split}
		\Phi(\bfs ; \bfth) = \bfw^\top N(\bfs;\bfth_N) + \frac{1}{2} \bfs^\top (\bfA^\top 
		\bfA)\bfs + \bfb^\top \bfs + c, \\  \quad \text{where} \quad\bfth = (\bfw, \bfth_N, \bfA, \bfb, c).
		\end{split}
	\end{equation}
	The inputs $\bfs{=}(\bfx,t) \in \R^{d+1}$ correspond to space-time, $N(\bfs;\bfth_N) \colon \R^{d+1} {\to} \R^m$ is an NN,
	and $\bfth$ contains trainable weights: $\bfw\,\,{\in}\,\,\R^m$, $\bfth_N\,\,{\in}\,\,\R^p$, $\bfA\,\,{\in}\,\,\R^{\gamma \times (d+1)}$, $\bfb\,\,{\in}\,\,\R^{d+1}$, $c\,{\in}\,\R$, where rank $\gamma{=}\min (10,d)$ limits the number of parameters in $\bfA^\top \bfA$. 
	Here, $\bfA$, $\bfb$, and $c$ model quadratic potentials, i.e., linear dynamics; $N$ models nonlinear dynamics.
	
	In our experiments, for $N$, we use a simple two-layer residual neural network (ResNet)~\cite{he2016deep}	
	\begin{equation} \label{eq:ResNet}
		  \begin{split}
		  \bfa_0 & = \sigma(\bfK_0 \bfs + \bfb_0) \\ 
		  N(\bfs;\bfth_N) & = \bfa_0 + \sigma(\bfK_1 \bfa_0 + \bfb_1),\\ 
		  \end{split}
	\end{equation}	
	for $\bfth_N{=}(\bfK_0,\bfK_1, \bfb_0,\bfb_1)$ where $\bfK_0 \in \R^{m \times (d+1)}$, $\bfK_1 \in \R^{m \times m}$, and $\bfb_0, \bfb_1 \in \R^{m}$.
	We use the element-wise nonlinearity $\sigma(\bfx){=}\log(\exp(\bfx) + \exp(-\bfx))$, the antiderivative of hyperbolic tangent, i.e., $\sigma'(\bfx){=}\tanh(\bfx)$~\cite{ruthotto2020machine,onken2020otflow}.

	We solve the ODE-constrained optimization problem~\eqref{eq:full_opt} using the discretize-then-optimize approach~\cite{gholami2019anode,onken2020do}, where we define a discretization of the ODE, then optimize on that discretization.
	The model's forward pass uses a Runge-Kutta 4 integrator with $n_t$ time steps to approximate the constraints~\eqref{eq:ODE_cons}. 
	The objective function is then computed, and automatic differentiation~\cite{nocedal2006numerical} calculates the gradient of the objective function with respect to $\bfth$. We use ADAM~\cite{kingma2014adam}, a gradient-based stochastic method with momentum, to update the parameters $\bfth$. We iterate this process a selected number of times. For the \emph{learning rate} (step size) provided to ADAM, we follow a piece-wise constant decay schedule, e.g., we divide the learning rate by 10 every 800 iterations (Fig.~\ref{fig:penalizers}).
	
\begin{figure*}[t]
\centering

\begin{subfigure}{0.245\linewidth}
\centering \captionsetup{width=0.85\linewidth}%
  \includegraphics[width=\linewidth]{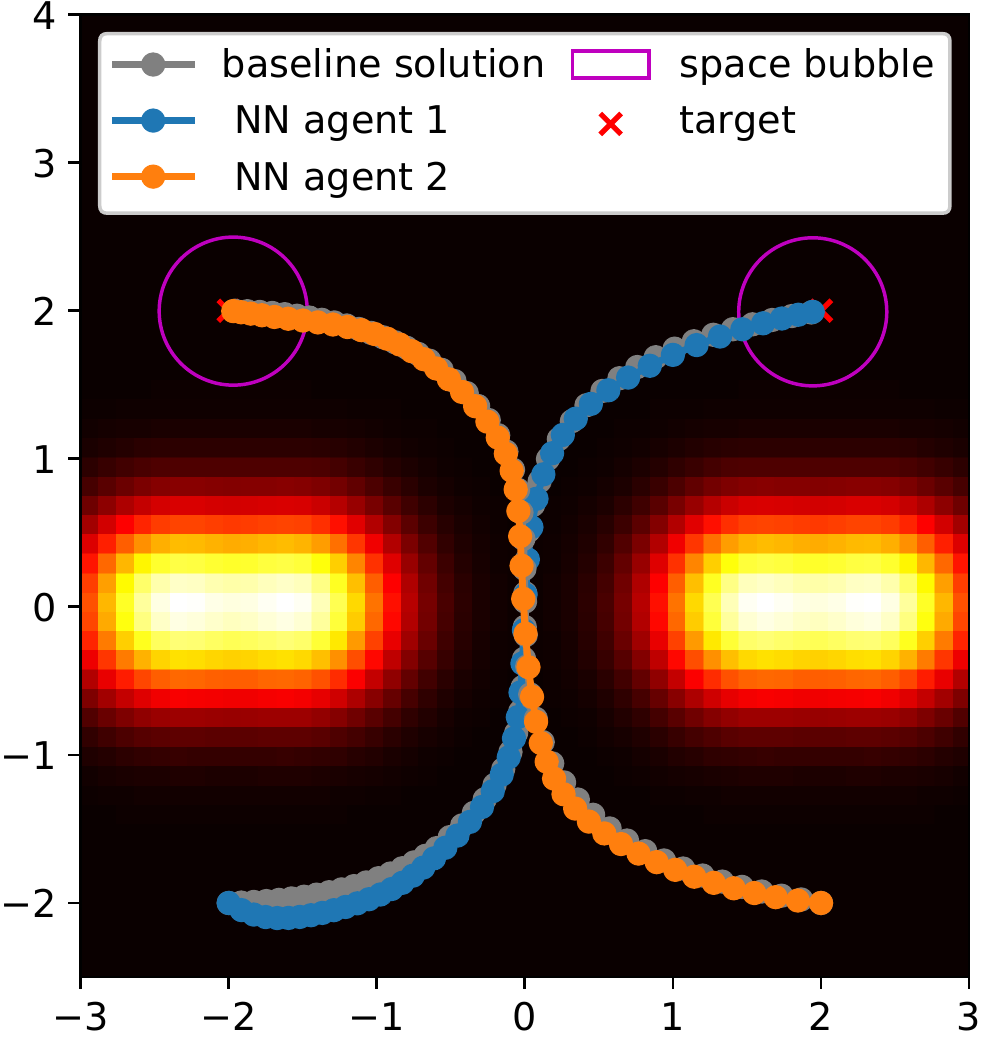}
	\subcaption{No shocks. The baseline and NN solutions for $\bfx_0$.}
	\label{fig:noshocks}
\end{subfigure}
\begin{subfigure}{0.24\linewidth}
  \centering \captionsetup{width=0.85\linewidth}%
   \includegraphics[width=\linewidth]{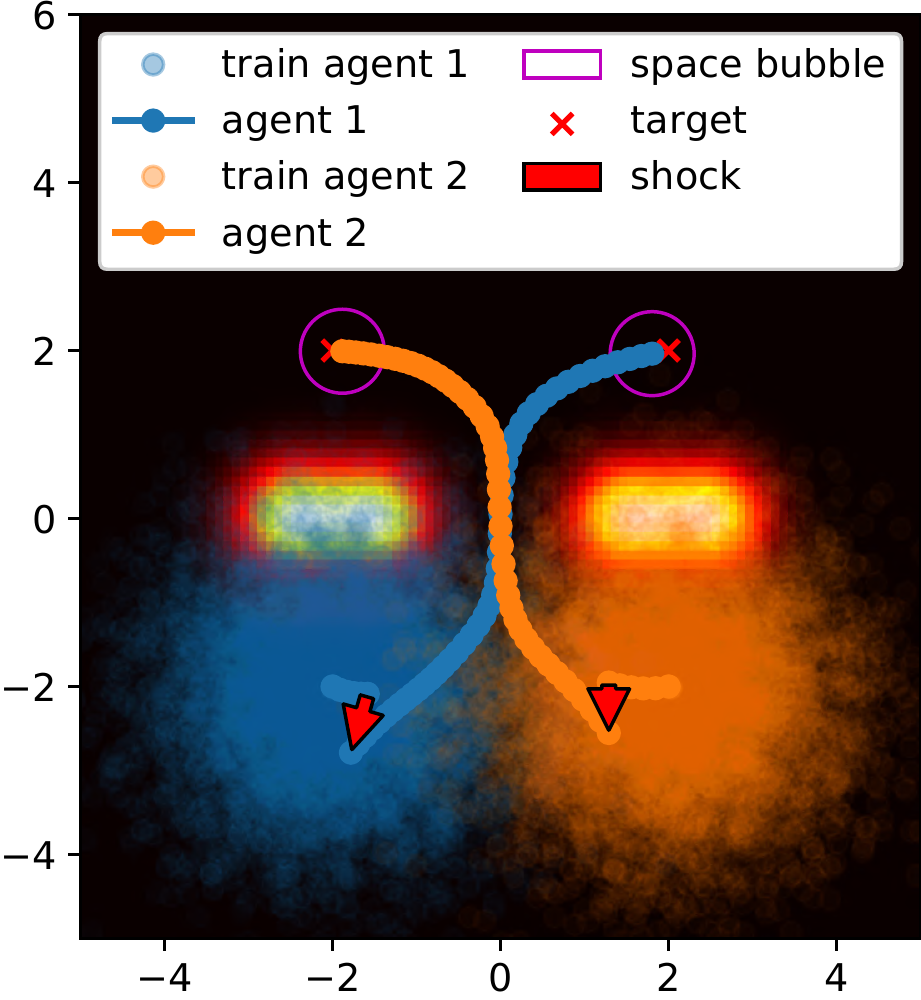}
  \subcaption{Minor shock $|\bfxi|{=}0.94$ within the training space.}
  \label{fig:minor_shock}
\end{subfigure}
\begin{subfigure}{0.24\linewidth}
  \centering \captionsetup{width=0.85\linewidth}%
   \includegraphics[width=\linewidth]{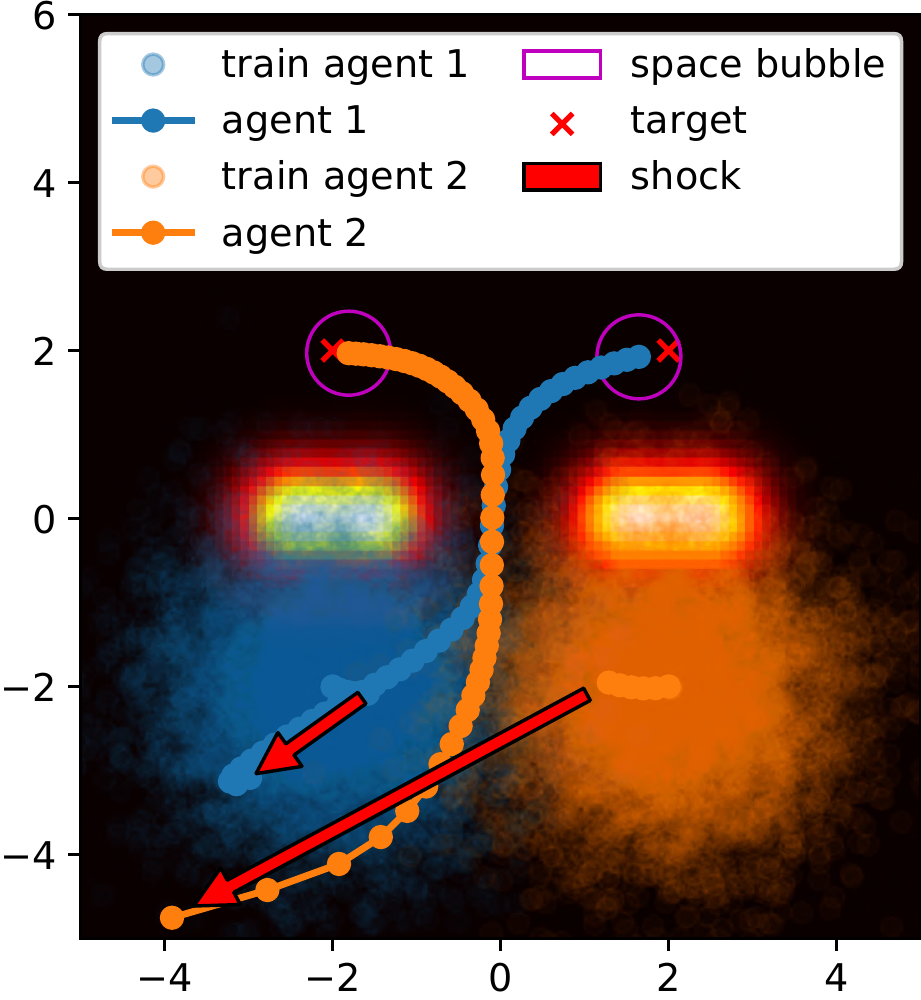}
  \subcaption{Major shock $|\bfxi|{=}6.2$ outside of the training space.}
  \label{fig:major_shock}
  \end{subfigure}
\begin{subfigure}{0.24\linewidth}
  \centering 
  \vspace{8pt}
  $\hat{\bfx} = \bfz_{\bfx_0}(0.1) + \bfxi$, $t= [0.1,1]$\\
  \vspace{4pt}
  \captionsetup{width=0.85\linewidth}%
	\includegraphics[width=\linewidth]{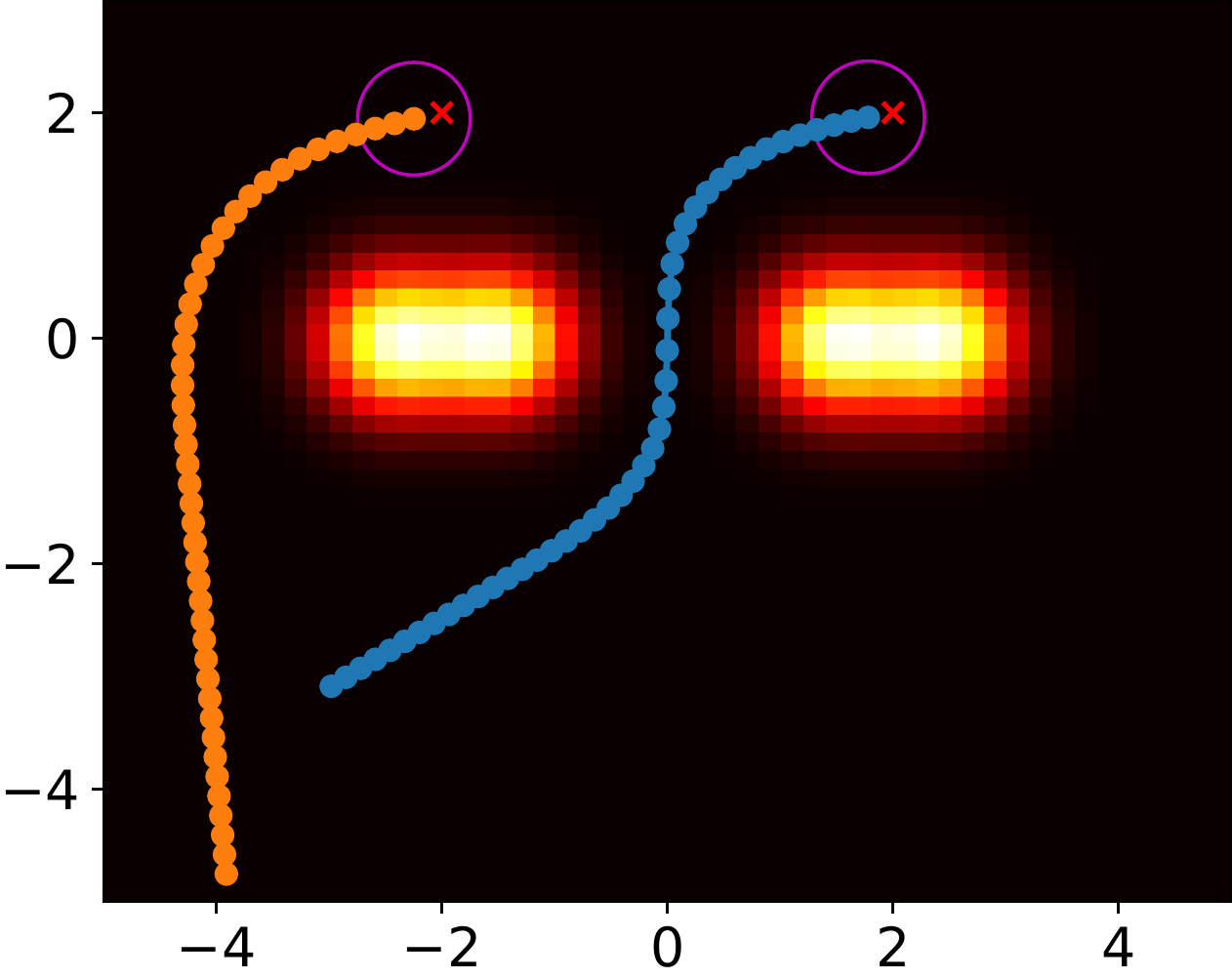}
  \subcaption{Baseline solution for (0.1,$\hat{\bfx}$) after major shock.}
  \label{fig:baseline_shock}
  \end{subfigure}  
\caption{Corridor problem. The NN handles a shock $\bfxi$ at time $t{=}0.1$ (depicted with red arrows) along the trajectory. Accompanying videos are available at \url{https://imgur.com/a/HWqlIot}.}
\label{fig:shock}
\end{figure*}

	The number of time steps $n_t$ is selected \textit{a priori} as a model hyperparameter. Large $n_t$ leads to high computation and training time while reducing error; meanwhile, too small $n_t$ leads to overfitting to a refinement of the time discretization of the trajectories. To avoid overfitting, we use more time steps for the hold-out validation set. For the corridor problem (Sec.~\ref{sec:softcorridor}), we use $n_t{=}20$ for training and $n_t{=}50$ for validation (Fig.~\ref{fig:noshocks}). For the swarm problem (Sec.~\ref{subsec:swarm}), we use $n_t{=}26$ for training and $n_t{=}80$ for validation (Fig.~\ref{fig:quadtraj}). Training on a single NVIDIA Quadro RTX 8000 GPU requires about ten minutes for the  corridor problem ($d{=}4$) and less than one hour for the swarm problem ($d{=}150$).
	
	Other hyperparameters include the width of the ResNet $m$ and the multipliers $\beta_1,\beta_2,\beta_3$. In contrast, some multipliers are inherent to the OC problem---e.g., $\alpha_1,\alpha_2,\alpha_3$---and are used for both the baseline and NN. Our Python implementation with all tuned hyperparameters is available at \url{https://github.com/donken/NeuralOC}.

\section{Numerical Experiments} 
\label{sec:numerical_examples}

We present a comparable baseline approach and solve two multi-agent OC problems of dimensions 4 and 150.

\subsection{Baseline: Discrete Optimization for a Single Initial State}

We provide a comparable local solution method that solves the OC problem for a fixed initial state  $\bfz^{(0)} = \bfx_0$.
Applying forward Euler to the state equation and a midpoint rule to the integrals, we obtain the discrete optimization problem
\begin{equation} \label{eq:baseline}
	\begin{split}
    \min_{ \{\bfu^{(k)}\} } \quad &
    \, G \left(\bfz^{(n_t)} \right) + h \sum_{k = 0}^{n_t-1}  L\left(s^{(k)},\bfz^{(k)},\bfu^{(k)}\right)\\
    \text{s.t.}  \quad  &\bfz^{(k+1)} = \bfz^{(k)} + h \, f\big(s^{(k)},\bfz^{(k)},\bfu^{(k)} \big),
    \end{split}
\end{equation}
where $h{=}T/n_t$. We use $T{=}1$ and $n_t{=}50$ and solve~\eqref{eq:baseline} using standard nonlinear programming techniques.

\subsection{Corridor Problem} \label{sec:softcorridor}

We design a four-dimensional problem in which two agents attempt to reach fixed targets on the other side of two hills. Suppose the agents with radius $r{=}0.5$ start at $x_{1}{=}[-2, -2]^{\top}$ and $x_{2}{=}[2,-2]^{\top}$ with targets $y_1{=}[2, 2]^{\top}$ and $y_2{=}[-2, 2]^{\top}$. Thus, the initial and target joint-states are
$\bfx_{0}{=}[-2,-2,2,-2]^{\top}$ and $\bfy{=}[2,2,-2,2]^{\top}$, respectively.
We sample from $\rho_0$, which is a Gaussian centered at $\bfx_0$ with an identity covariance. These sampled initial positions form the training set $\bfX$. We sample again to create the validation set.

The obstacles are defined by the spatio-temporal cost function $Q_i$, which we define in this experiment as the sum of four Gaussians.
The energy terms are given by
\begin{equation} \label{eq:Ei}
	E_i \big( u_{x_i} \big) = \hf \|u_{x_i}\|^2, 
\end{equation}
and the dynamics are given by $f(t,\bfx,\bfu) = \bfu$; the controls are the velocities.
We model interactions via
\begin{equation} \label{eq:Wij}
	W_{ij}(z_{x_i},z_{x_j}) =
	\begin{cases}
     \exp \left( -  \frac{ \| z_{x_i} - z_{x_j} \|_2^2}{2r^2} \right) ,  &\| z_{x_i} {-} z_{x_j} \|_2 < 2r, \\
     0 , &\text{otherwise}.
	\end{cases}
\end{equation}
for agents with radius $r$.
For terminal costs, we choose
\begin{equation} \label{eq:G}
	G\big(\bfz_{\bfx}(T)\big) = \frac{\alpha_1}{2}\| \bfz_{\bfx}(T) - \bfy \|^2.
\end{equation}
We select multipliers $\alpha_1{=}100$, $\alpha_2{=}10,000$, and $\alpha_3{=}300$.

The baseline and the NN solve the problem with comparable trajectories (Fig.~\ref{fig:noshocks}) and accuracy (Table~\ref{tab:softcorridor}). 
The NN achieves a marginally worse $\ell$ value but slightly better $G$  value (Table~\ref{tab:softcorridor}). We attribute the $G$ improvement to the final-time HJB penalizers (Fig.~\ref{fig:penalizers}). 

\subsection{Shocks}\label{subsec:shocks}

We observe that approximating the value function leads to a shock-robust model. Since the controls are in the feedback form $\partial_s \bfz=-\nabla_{\bfp} H(s,\bfz,\nabla \Phi(s,\bfz))$, the model can quickly calculate updated trajectories despite shocks to the system.
As an example, we consider a shock $\bfxi \in \R^d$ (implemented as a random shift) to the system at time $t{=}0.1$ when solving the corridor problem for $t\in[0,T]$ (Fig.~\ref{fig:shock}).

Our method is designed to handle minor shocks that stay within the space of trajectories of the initial cloud about $\bfx_0$. Our model computes a trajectory to $\bfy$ for many initial points. Therefore, for point $\widetilde{\bfx} \in \bfX$, the model provides dynamics $f(t,\bfz_{\widetilde{\bfx}}(t), \bfu_{\widetilde{\bfx}}(t))$ before the shock. After the shock, the state picks up and follows the trajectory of some other point $\hat{\bfx} \in \bfX$ (Fig.~\ref{fig:minor_shock}). Thus, the total trajectory has two portions
\begin{equation*} \label{eq:shock}
	\begin{split}
		&\bfz_{\widetilde{\bfx}}(0.1) = \int_0^{0.1} f \big( t,\bfz_{\widetilde{\bfx}}(t),\bfu_{\widetilde{\bfx}}(t) \big) \, \du t, \quad \bfz_{\widetilde{\bfx}}(0)= \widetilde{\bfx}, \quad \text{and} \\
		&\bfz_{\hat{\bfx}}(1) = \int_{0.1}^{1} f \big( t,\bfz_{\hat{\bfx}}(t),\bfu_{\hat{\bfx}}(t) \big) \, \du t, \, \bfz_{\hat{\bfx}}(0.1)=\bfz_{\widetilde{\bfx}}(0.1)+\bfxi ,
	\end{split}
\end{equation*} 
before and after the shock, respectively.
A minor shock can thus be categorized as switching from one trajectory to another. (Fig.~\ref{fig:minor_shock}). The NN and baseline results of the control problem along $t{=}[0.1,1]$ are similar (Table~\ref{tab:softcorridor}).

\begin{table}[t]
  \centering
  \caption{Comparison for single instance $\bfx_0$.}
	\begin{tabular}{clccc}
	   	\toprule
	 	 Scenario & Method & $\ell + G$ & $\ell$ & $G$ \\
	   	\midrule
	   	no shocks & NN       & \hphantom{1}62.19 & \hphantom{1}61.98 & 0.21 \\
	   	$t \in [0,1]$ & Baseline & \hphantom{1}61.33 & \hphantom{1}61.02 & 0.31 \\
	     \midrule
	     after shock $|\bfxi|=0.94$ & NN  & \hphantom{1}60.54 &  \hphantom{1}60.34  &  0.20  \\
	     $t \in [0.1,1]$ & Baseline & \hphantom{1}59.79 & \hphantom{1}59.46 & 0.33  \\ 
	     \midrule 
	     after shock $|\bfxi|=6.2$ & NN       & 151.67 &  150.63 &  1.03\\
	     $t \in [0.1,1]$ & Baseline & \hphantom{1}71.77 & \hphantom{1}71.22 & 0.55 \\
	 	\bottomrule
	\end{tabular} 
	\label{tab:softcorridor}
\end{table}

Interestingly, our model extends outside the training region (Fig.~\ref{fig:major_shock}).
Although the vast majority of NNs cannot extrapolate, our NN still solves the control problem after a major shock, demonstrating some extrapolation capabilities. We note that the NN solves the original problem for $\bfx_0$ to near optimality. However, after such a large shock, the NN solves the control problem, but sub-optimally. In our example, we compare the NN's solution (Fig.~\ref{fig:major_shock}) with the baseline's solution for $t{=}[0.1,1]$ (Fig.~\ref{fig:baseline_shock}). The NN learned a solution where agent 2 passes through the corridor before agent 1. After the major shock, the NN applies these dynamics (Fig.~\ref{fig:major_shock}), while the baseline finds a more optimal solution (Fig.~\ref{fig:baseline_shock}). The NN is roughly 100\% suboptimal (Table~\ref{tab:softcorridor}). 

We attribute the shock robustness to the NN parameterization of the global value function.
Experimentally, the shock robustness of our model (Fig.~\ref{fig:shock}) does not noticeably differ from a model trained without penalization.
Since the NN is trained prior and offline, it handles shocks in real-time. In contrast, methods that solve for a single trajectory---e.g., the baseline---must pause to recompute following a shock.

\subsection{Swarm Trajectory Planning} \label{subsec:swarm}

\begin{figure}
	\centering
	\includegraphics[width=0.95\columnwidth]{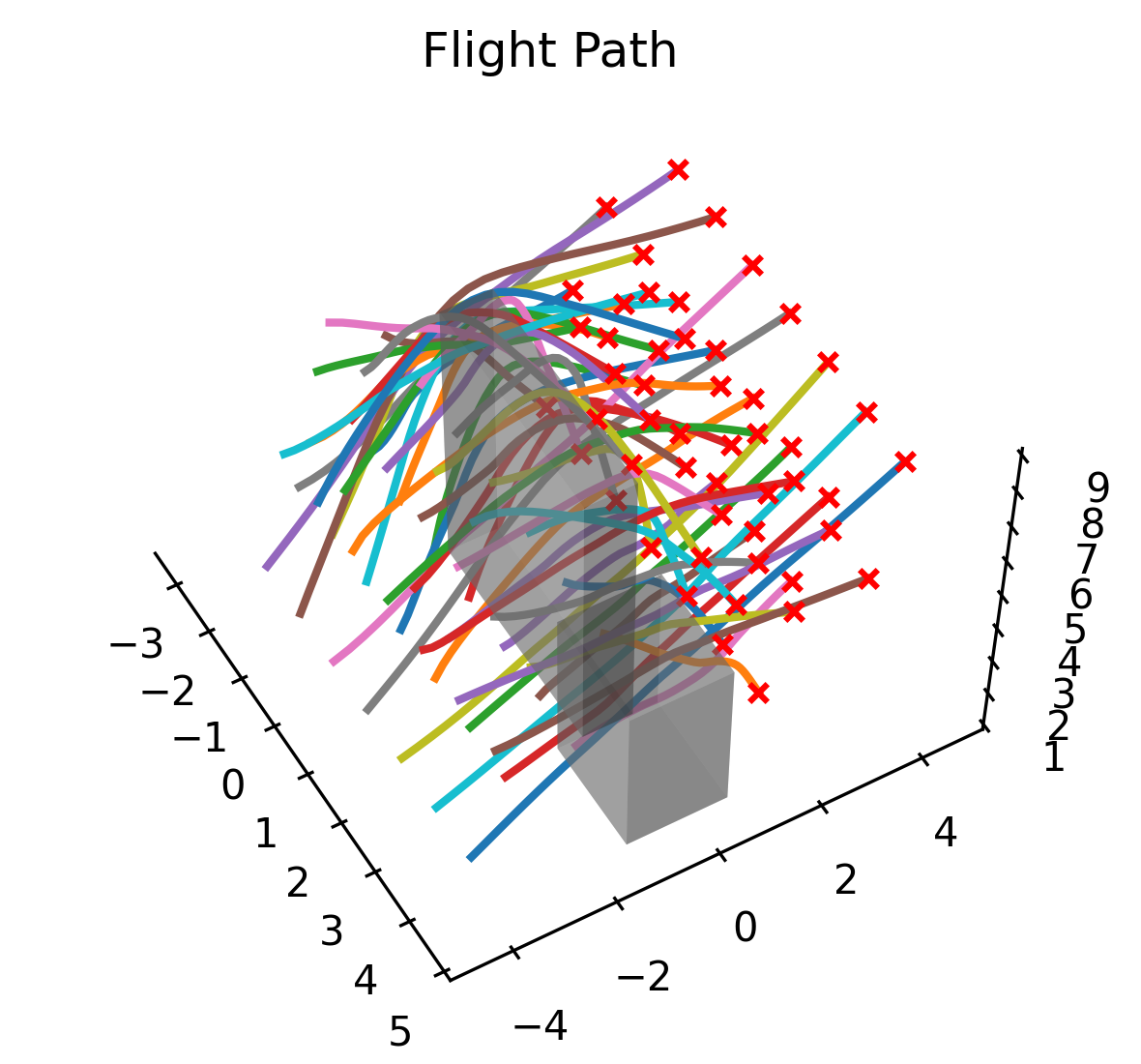}
	\caption{Swarm Trajectory Planning for 50 agents in $\R^3$.}
	\label{fig:quadtraj}
\end{figure}

We demonstrate the high-dimensional capabilities of our model by solving a swarm trajectory planning problem in the spirit of~\cite{honig2018trajectory}.
The swarm problem contains 50 three-dimensional agents that fly from initial to target positions while avoiding each other and obstacles. We construct $Q_i$ to model rectangular prism obstacles and use \eqref{eq:Ei}, \eqref{eq:Wij}, and \eqref{eq:G} for energy, interaction, and terminal costs. The NN guides all agents around the obstacles (Fig.~\ref{fig:quadtraj}). 
Naturally, if the time discretization is too coarse (small $n_t$), the model may simulate collisions solely due to inaccurate integration. 
In validation, we see that the agents avoid the obstacles and each other by observing values for $Q$ and $W$ are exactly 0.

\section{Conclusion and Outlook}

We formulate and demonstrate an efficient NN approach for solving high-dimensional OC problems. 
Our method aims at computing the optimal control in feedback form in the relevant subset of the space-time domain.
It combines the high-dimensional scalability from PMP and the global nature from HJB approaches.
Using a numerical example, we demonstrate that the obtained feedback form generalizes outside the training space, which allows the agents to react to unforeseen events such as shocks. Our future endeavors relate to further experimentation of our method on OC problems with more involved dynamics and blending our method (trained prior) with distributed approaches in deployment.

\addtolength{\textheight}{-12cm}   

\section*{Acknowledgment}
We thank Reza Karimi for assisting with figure creation.

\bibliographystyle{IEEEtran}
\bibliography{main}

\end{document}

%% file: abbrv.tex
\newcommand{\hf}{{\frac 12}}

\newcommand{\bfa}{\boldsymbol{a}}
\newcommand{\bfb}{\boldsymbol{b}}

\newcommand{\bfp}{\boldsymbol{p}}

\newcommand{\bfs}{\boldsymbol{s}}

\newcommand{\bfu}{\boldsymbol{u}}

\newcommand{\bfw}{\boldsymbol{w}}
\newcommand{\bfx}{\boldsymbol{x}}
\newcommand{\bfy}{\boldsymbol{y}}
\newcommand{\bfz}{\boldsymbol{z}}

\newcommand{\bfA}{\boldsymbol{A}}

\newcommand{\bfK}{\boldsymbol{K}}

\newcommand{\bfX}{\boldsymbol{X}}

\newcommand{\CP}{{\cal P}}

\newcommand{\bfxi}{{\boldsymbol \xi}}

\newcommand{\R}{\ensuremath{\mathds{R}}}
\newcommand{\E}{\ensuremath{\mathds{E}}}

\def\du{\ensuremath{\mathrm{d}}}

\newcommand{\bfth}{\boldsymbol{\theta}}